\documentclass[11pt]{article} 

\usepackage{latexsym}
\usepackage{amssymb}
\usepackage[mathscr]{eucal}

\begin{document} 

\begin{center} 

 {\Large\bf     G\"{o}del, Tarski, Church, and the Liar} 

\vspace{3mm}

                      {\bf    Gy\"{o}rgy Ser\'{e}ny} 

\end{center} 

\vspace{7mm}

\hspace{2mm}   `What is a sadist? A sadist is a person who is kind to a 
masochist.' When I first heard this joke, I smiled, of course, 
realizing that the {\it brutality} of the sadist is, in fact, 
{\it kindness} to 
a masochist, which is the source of the comic effect being 
a{\small (}$\!$n apparent$\!${\small )} 
contradiction resolved. What I did not realize was that I would have 
had another reason to smile: the {\it kindness} of the sadist 
{\it tortures} the 
masochist. Well, then how should a sadist treat a masochist? Is this 
really a joke or might it be a serious scientific question?
\footnote{Frigyes Karinthy, whose name in Hungary means something like 
  those of Oscar Wilde and Monty Python together in Britain, posed the 
  following question: 

\vspace{0.3mm}

{\scriptsize 
  If to be a madman is to have an obsession, and my only obsession 
  is that I am a madman, have I really lost my mind or not?}} 
The writer Arthur Koestler in his study, examining the connection 
between science and art, thinks that it is both ([K] p.95.):

\vspace{1mm} 

{\small Comic discovery is 
paradox stated -- scientific discovery is paradox resolved.} 

\vspace{1mm}

I think, he is utterly right. In fact, it will emerge from what we 
shall do below that the general ideas underlying some of the main 
mathematical results of this century, those concerning the 
incompleteness and undecidability of arithmetic and the undefinability 
of truth within it can, at least partly, 
be taken as different ways to resolve the archetype of our initial 
paradoxical question. The paradoxical nature of the famous 
G\"{o}del's theorem\footnote{Using 
W.V.Quine's words: 

\vspace{0.3mm}

{\scriptsize [G\"{o}del's discovery] is decidedly paradoxical in the sense that 
it upsets crucial preconceptions. We used to think that mathematical 
truth consisted in provability. ([Q] p.17.)}}
and its less broadly known but equally important and enlightening 
relatives, the theorems of Tarski and Church, which together 
constitute the core of the family of the results of modern 
mathematical logic describing the theoretical limitations of formal 
reasoning, is inherited from their ancient and not less famous 
ancestor, one of the most important of all logical paradoxes, that of 
the Liar.
\footnote{The fact that G\"{o}del's theorem and the Liar paradox are 
  related closely is, of course, not only well known, but is a part 
  of the common knowledge of logician community. Actually, almost every 
  more or less formal treatment of the theorem makes a reference to 
  this connection. Well, G\"{o}del himself is not an exception as he made a 
  remark in the paper announcing his result (cf. [G]): 

\vspace{0.3mm}

  {\scriptsize The analogy between this result and Richard's antinomy leaps to the 
  eye; there is also a close relationship with the `liar' antinomy, 
  since ... we are ... confronted with a proposition which asserts its 
  own unprovability.} 

\vspace{0.3mm}

  In the light of the fact that the existence of this connection is a 
  commonplace, all the more surprising that very little can be 
  learnt about its exact nature except perhaps 
  that it is some kind of similarity or analogy.}

\hspace{3mm}In fact, our aim in this article is just to show 
that an abstract 
formal variant of the Liar paradox constitutes a general 
conceptual schema that, revealing their common 
logical roots, connects the theorems 
referred to above and, at the same time, demonstrates 
that, in a sense, 
these are the only possible 
relevant limitation theorems formulated in terms of truth and 
provability alone that can be considered 
as different manifestations of the Liar paradox. On the other hand, 
as we shall illustrate by a simple example, this abstract version 
of the paradox 
 opens up the 
possibility to formulate related results 
concerning notions other than just those of the 
truth and provability. 
  
\hspace{3mm}To obtain 
the resolution of the paradox in a general, purely formal wording, 
we shall reformulate 
the infamous statement of the Liar in a step by step manner 
in two main stages.
First we shall seek an ordinary--language equivalent of the paradox 
in a form that shows clearly its logical structure, then we shall 
directly translate the expression we have obtained into 
a formal language. By  applying   
the generalized version of 
the result of this formalization 
process 
to the logical systems used to describe mathematical 
theories, we shall finally  
 make explicit the connection between the age old paradox and 
its modern reincarnations. Let us start with the origins.

\pagebreak

\vspace{1mm}

\begin{center} 
 {\large\bf                                1. The paradoxes  }

\end{center}

{\bf a. The Liar paradox}

\vspace{1mm}

  \hspace{2mm} Perhaps the 
 most widely known version of the paradox, which is called 
the paradox of Epimenides, appears in the New Testament. St.Paul, 
referring to Epimenides of Crete, says \linebreak (St. Paul's epistle to Titus 
(I, 12)): 

\vspace{1mm}

{\small One of themselves, even a prophet of their own, said, `Cretans 
are always liars, evil beasts, lazy gluttons.' }

\vspace{1mm}

As a matter of fact, the statement `Cretans are always liars' uttered by 
a Cretan is not a logical paradox since it must simply be false. On 
the other hand, it indeed exhibits remarkable paradoxical features since 
the apparently contingent fact that some Cretans existed who sometime 
told the truth turns out to be a logical necessity.

\hspace{2mm}  The paradox in a clear-cut purely logical form is generally 
attributed to the Greek philosopher Eubulides from Miletos. He lived 
in the fourth century B.C. and belonged to the school of Megara. 
Eubulides formulated the paradox in the form of the following puzzle\,:
\footnote{The fact that the Megarians did take puzzles seriously is 
  reflected in the story according to which one of them, Diodorus Cronus 
  allegedly committed suicide because he was not immediately able to 
  solve a logical puzzle.}

\vspace{1mm}

   {\small A man says that he is lying. Is what he says true or 
false?}

\vspace{1mm}

   So the paradox consists in the fact that the 
sentence `I am lying' can be neither true nor false since its truth implies its 
falsity and {\it vice versa}.       
   Now, in logical investigations it is natural to consider the 
impersonal version of the original paradoxical sentence. Consequently, 
from now on, we shall examine the sentence below, which will simply be 
called the {\it Liar}\,:

\vspace{1mm}

  \hfill         This sentence is false.   \hfill \

\vspace{1mm}

\hspace{2mm}   One of the innumerable attempts to get rid of the paradox 
caused by 
the impossibility to assign truth value to the Liar was the claim according 
to which it only pretends to be a declarative sentence of the 
subject--predicative form, which must be either true or false, since it 
has not really got a subject. In fact, in order for a noun phrase to 
play the role of the subject of a meaningful sentence, it has to be a 
complete linguistic entity. In this case, however, the noun phrase 
under consideration, the sentence itself, is just {\it in statu 
nascendi}
\footnote{being born} at the point where the need for the subject 
appears. Thus, without a proper subject, the sentence is not 
paradoxical but simply meaningless. Alternatively, we can eliminate 
the paradox accepting the analysis of Quine:  

\vspace{1mm}

{\small \ldots the phrase `this sentence', so used refers to nothing. This is 
claimed on the grounds that you cannot get rid of the phrase by 
supplying a sentence that is referred to. For what sentence does the 
phrase refer to? The sentence `This sentence is false'. If, 
accordingly, we supplant the phrase `This sentence' by a quotation of 
the sentence referred to, we get: `\,``This sentence is false'' 
is false'. 
But the whole outside sentence here attributes falsity no longer to 
itself but merely to something other than itself, thereby engendering 
no paradox. ([Q] p.7)}.
\footnote{A similar criticism is propounded by G. Ryle (cf. [R]): 

\vspace{0.3mm}

{\scriptsize \noindent  If unpacked our pretended assertion would run 
`The current statement 
  \{namely that the current statement [namely that the current  
  statement (namely that the current statement \ldots '. The brackets  are 
  never closed; no verb is ever reached; no statement of which we  can 
  even ask whether it is true or false is ever adduced. }

\noindent  I think that the infinite regression referred to in this quotation 
  can be represented in the physical world, of course, by the   very nature of the representation, only in a non--complete way. Let us try 
  to represent a picture that is about itself. Well, the optical phenomenon 
  it yields might be like the infinite--looking array of ever 
  decreasing pictures the sight of which we experience when, standing 
  between two mirrors facing each other, we are looking at one of 
  them.}

\vspace{1mm}

Anyway, the Liar is in this form untenable, so we should modify it in order to 
preserve our paradox.$\!$\footnote{Obviously, every variant of the usual practice of giving a name 
  to our sentence and referring to it by this name, that is e.g. the 
  following presentation of the paradox

\vspace{0.3mm}

  (*)  \hfill    The sentence (*) is false , \hfill \

\vspace{0.3mm}

  does not make much difference and is open to essentially the same   
  objections as its original version.}
Since, clearly, the paradoxical character of the Liar 
stems from its attributing falsehood to {\it itself}, we can 
get help from the 
analysis of another famous paradox explicitly based on the self--reference
\nolinebreak
\footnote{ 
the linguistic phenomenon of being about itself} providing 
opportunity for a more systematic examination of self--referential 
sentences.\footnote{As we shall see 
later, however, there are both informal and formal arguments proving 
that, in itself, the self-reference in the Liar is innocent; it 
cannot be 
blamed for the paradox. Actually, our formalization process will 
clarify the difference between the roles played by different factors 
the paradox is built on.}

\vspace{2mm}

{\bf b. The Quine paradox  }

\vspace{1mm}

\hspace{2mm}   
The {\it paradox of heterologicality} was devised by the German 
mathematician Kurt Grelling in 1908. An adjective is called 
{\it autological} if the property denoted by the adjective holds for the 
adjective itself, otherwise it is called {\it heterological}. For example, 
`short', `English', `polysyllabic', `mispelt', `adjectival' are 
autological, while `long', `German', `monosyllabic', `obsolete', 
`obscene', and `non--adjectival', for that matter, are heterological. 
Well, then, what about the adjective `heterological' itself? Clearly, 
it is neither true nor false of itself so it can be neither 
autological nor heterological, providing a paradox.

\hspace{2mm}  The Grelling paradox 
is closely related to the Liar. 
Actually, the 
fact that `heterological' is neither 
true nor false of itself can be 
utilized to formulate a sentence with 
the same paradoxical nature as 
the Liar. Indeed, for any adjective, say `white', `consists of four words', 
or `heterological', the fact that an 
object possesses the property 
referred to by the adjective can be expressed by a sentence 
resulting from the application of 
the adjective to the noun naming 
the object. This sentence is 
true just in the case that the object does possess 
the property considered. 
Thus `snow is white' is true, `\,``snow is white'' consists of four 
words ' and `\,``does not consist of four words'' 
consists of four words' are 
false, while 

\vspace{1mm}

(1) \hfill  `consists of four words' consists of four words  \hfill \

\vspace{1mm}

is again true.\footnote{The existence of completely 
meaningful and immaculate 
applications of autologicality like this sentence provides 
informal evidence that the self-reference in the Liar, in itself, 
cannot be the source of the paradox. Nevertheless, self-reference is 
not so innocent and  banal a  phenomenon as, for example, Hofstadter
 \mbox{thinks it is ([H] p.7):} 

  {\scriptsize Self-reference is ubiquitous. It happens every time anyone says `I' 
  or `me' or `word' or `speak' or `mouth'. It happens every time a 
  newspaper prints a story about reporters, every time someone writes 
  a book about writing, designs a book about book design, makes a 
  movie about movies, or writes an article about self-reference.} 

  Well, the examples listed in this quotation, with the only exception 
  of the personal pronoun, are {\it not} self-referential (at least in the   
  sense of the word related to interesting logico-linguistic   
  phenomena worth considering), since they are simply not about 
  themselves, the word `word' refers to any word not this particular 
  one (`these words' -- but {\it not} `this word' -- 
 would do something like 
   that), a story about reporters is 
  not a story about the particular story concerned, a book about 
  writing is not a book about itself (but about writing) etc. On the 
  other hand,  clearly, proper self-reference is not a boring 
  everyday linguistic act  as the famous paradoxes (among them 
  those we are just examining) exemplify, which are all built on it in 
  some way or other.} So far so good. But what about the following 
  sentence\,:

\vspace{1mm}

(2)   \hfill         `heterological' is heterological.  \hfill \

\vspace{1mm}

Clearly, it can be neither true nor false, yielding a variant of the 
Grelling paradox in a form of a sentence sharing the main 
characteristic of the Liar being without truth value.
 Yet, it is immune against the kind of criticism that 
the Liar was open to since it does have a proper subject. Of course, it 
relies heavily on the meaning of the technical term in it, so, in 
order to obtain a sentence that stands on its own feet, the term 
`heterological' must be eliminated. Before looking for a 
common language version of (2), however, we should step back and 
examine the autological sentence (1) more thoroughly. For the way it 
expresses its meaning points to an underlying general pattern 
for handling 
self-reference in an apparently unobjectionable way. In fact, it contains implicitly the following definition of self-reference\,:

\vspace{1mm}

(3) \hfill  `applied to itself' \  is equivalent to \ 
            `appended to its own quotation'.   \hfill \

\vspace{1mm}

Now, by definition,  
  `heterological'   means    
 `yields falsehood when applied to itself', which, in turn, by (3),    
 is equivalent to     
     `yields falsehood when appended to its own quotation'. 
Applying this translation to (2), we obtain a masterpiece 
version of the Liar due to Quine \mbox{(cf. [Q] p.7).} 

\vspace{1mm}

The {\it Quine paradox}\,:

\vspace{1mm}

\ \hfill \begin{minipage}[t]{10cm}
        `yields falsehood when appended to its own quotation' 

          \hspace{1mm}yields falsehood when appended to its own quotation. 

\end{minipage} \hfill \

\vspace{1.5mm}

In fact, this sentence does have a subject, which is just the name of the string of words in its second row. For, according to the usual 
practice (which we also tacitly assumed already),  
the name of a string of words is simply the string itself between 
quotation marks. Thus the Quine paradox tells something about the 
bearer of this name, that is, the string of words in its second row, 
namely that the sentence that can 
be obtained by putting it down in a particular way -- first between quotation marks then 
without them -- is false. But, of course, the resulting sentence is the 
Quine paradox itself. That is, this sentence indeed says of itself 
that it is false. Moreover, having a proper subject, it is a 
completely well-formed sentence. 

\hspace{2mm}   Having obtained the Quine paradox, it seems that we 
achieved our 
aim to present a natural language version of the Liar in a form that 
can be formalized in a straightforward way. As a matter of fact, we indeed made an 
important step in the right direction, but, as we shall see below, 
there remains yet another step to take since the paradox above 
 has  
some essential imperfections. The way 
  we shall improve its formulation is implicitly present in 
 the critical arguments this claim is based on. Actually, it is a 
reassuring fact that they all point in the same direction.

\vspace{5mm}

{\bf c. The Findlay paradox}

\vspace{1mm}

\hspace{1mm}
First of all, having given due credit to the ingenuity of the Quine 
paradox, we should notice that we have been cheated by it. There is 
a trick that made it possible to formulate it as a well--formed 
sentence without using the notion of heterologicality. On the one 
hand, it is clear that the bearers of truth and falsity are the 
declarative sentences; they are those linguistic entities that can 
be true or false.
\footnote{An important remark is in order here. We shall not 
  make the usual distinction made by philosophers of language between 
  the three notions that can appear in connection with this question, 
  namely the notion of a sentence, as the linguistic embodiment of its 
  meaning, that of the meaning itself, which is called a proposition 
  (expressed by the sentence), and finally the notion of a statement, 
  which is a word usually used to refer to the act that we perform in 
  uttering a declarative sentence. The reason is that, although generally 
  these distinctions indeed may be relevant (cf. e.g. [A]), from the 
  point of view of our analysis they are entirely immaterial.}
This means that 
the Quine paradox can be reformulated as follows: 

\vspace{1mm}

(4) \hfill \begin{minipage}[t]{12cm}  
`yields a false statement when appended to its own quotation' 

\hspace{1mm}yields a false statement when appended to its own 
quotation.

\end{minipage} \hfill  \

\vspace{1mm}

The subject of this sentence, the actor who performs the action 
expressed by the verb `yield' is the following string of words 

\vspace{1mm}

(5) \hfill  yields a false statement when appended to its own quotation.
    \hfill \ \ \ \ \ \ \ \ \ \  \ \ \ 

\vspace{1mm}

But it is an abuse of the verb `yield' since, clearly, it is not this 
string of words that makes the false statement! The reason is simple 
enough: there is no string of words that can build a sentence, true or 
false. Any sentence can come to existence only if there is an actor 
(supposedly intelligent), perhaps unknown, who performs the act of 
building it the way described in the Quine paradox, that is, 
 attaching the string of words considered to its quotation. 
What is more, it is highly probable that this  criticism can be 
generalized to any possible reformulation of (5), since it 
seems self--evident that there is no linguistic means by which 
 a part 
of a sentence as a linguistic entity can render the whole sentence 
whose part it is, false (or true, for that matter). The reason is the 
same as above, the sentences are not {\it rendered} false or true 
by their parts, 
they {\it are} false or true.
\footnote{Similar arguments apply to analogous but less neat 
  reformulations of (2) as e.g. the following one 
\mbox{(cf. [R] p.167):}

\vspace{0.3mm}

       \noindent  {\scriptsize `appended to its own quotation is false' 
           appended to its own quotation is false, }

\vspace{0.3mm}

 \noindent which makes the deception more apparent since, 
 obviously, it states 
  that the {\it phrase} `appended to its own quotation is false' is false 
  (when appended to its own quotation), which is plainly absurd since 
  this phrase is neither false nor true, it is simply meaningless in 
  itself. One can hardly suppress the feeling that the language cannot 
  be deceived: what cannot be said, cannot be thought. Well, 
what cannot be said in a natural language may be formulated 
in an artificial one.  For within a formalized system, the 
  linguistic competence (the ability to understand and use everyday 
  language correctly) does not matter. In fact, interpreting the 
  phrase `applied to itself' as meaning `followed by its own 
  quotation' and mirroring formally the latter by concatenation, 
  Smullyan has managed to construct formal languages capable of  
  handling self-reference (see [S1]). Concatenation was, among others, 
  used by Quine as a basis for various formal constructions (cf. [Q1]) 
  and it is probable that his paradox formulated in 
  everyday English was influenced by his artificial language practice. 
  On the other hand, the idea underlying Smullyan's formal 
  construction is related to Quine's informal one (cf. [S1] p.56 
  footnote 6). Formal systems, however, endowed  by 
     concatenation are not widely used. I wonder whether 
  the reason why the substitution, 
  as a formal linguistic tool 
  (which we shall describe soon), is preferred to   
  concatenation is 
  the fact that substitution  (as opposed to concatenation)   is commonly used in 
  mathematics, and, 
  clearly, the best way for a  formal language to satisfy the 
  natural requirement to be as easily comprehensible as possible is 
  to mimic formally the everyday languages, among them the language 
  of everyday mathematics, as faithfully as 
  possible.}
The moral is that we should try to construct another statement analogous 
to the Quine paradox using a {\it whole sentence} instead of the 
adjectival phrase (5).

\hspace{2mm}   Secondly, we may not expect that, without any modification, 
the Grelling paradox, which is a statement about an adjective, will 
provide a correct formulation of the Liar paradox, which is a 
statement about a sentence. A natural way to achieve this goal is, 
obviously, the same as we indicated above, that is, to formulate a 
 sentence that will play the role of the adjectival 
\mbox{phrase (5).} 

\hspace{2mm}  Finally, let us consider the most important 
 reason 
to change the 
phrase (5) in  
(4) to a suitable {\it sentence}. In fact, bearing in mind that our aim is to 
formalize the Liar, we should not have been satisfied with the Quine paradox  
even if the arguments above had not been convincing enough, since the 
fact it conveys the intended meaning is based on a contingent 
characteristic of English language. In fact, in English language, the 
usual word--order in declarative sentences is `subject + verb + \ldots' , 
that is, the role of a phrase in a sentence is shown not by its 
form but by its position in the sentence. Consequently, generally, if 
there is no indication to do otherwise, a noun--phrase at the beginning 
of a sentence, however complicated it may be, is automatically taken 
as the subject of the sentence by a speaker of English. Thus, in 
English, the linguistic application of an adjective (or adjectival 
phrase) to an object means a simple attachment of the adjective (or 
adjectival phrase) to the name of the object. (In the case of the 
Quine paradox, the name of the object is the object itself between 
quotation marks.) Therefore, the construction underlining the Quine paradox is based on 
a special feature of English language.\footnote
{I cannot resist the temptation to make use of my knowledge 
of an `exotic' language in order to give evidence backing this claim. 
In fact, the construction leading to the Quine paradox cannot 
be executed in Hungarian since a phrase like the Quine paradox is 
simply ungrammatical (though comprehensible) without a definite 
article preceding the subject of the sentence, that is, (3) is simply 
not true in Hungarian.} This fact means that we cannot expect 
the Quine paradox to exhibit its language--independent  
 logical structure, which is what we actually need 
in order to formulate our original paradox in a purely formal system. 
  Consequently, we ought to look for an abstract `language--free' 
procedure 
to execute the linguistic application of an adjective 
to a noun. 

\hspace{2mm}Fortunately, all languages as well as every logical system have a 
common formal tool for handling the application of one linguistic 
entity (an adjective (or adjectival phrase)) to another one (a noun 
(phrase)), namely the {\it substitution.} Given a sentence asserting that 
the indefinite subject of the sentence, denoted usually by a single 
letter and called (logical) variable, possess the property expressed 
by the adjective(al phrase) we would like to apply to our object, we 
can execute this linguistic application by substituting the name of 
our object for the indefinite subject of the sentence. 
 Thus e.g.  
 the 
linguistic application of the adjective `red' to 
snow can be executed by substituting the name of the latter 
for the letter $x$ in the sentence `$x$ is red', which will yield the sentence 
`snow is red'. Actually, we can alternatively interpret (1) as a result 
of a suitable substitution rather than that of appending a 
phrase to its name.\footnote{Note that not only the artificial languages use symbols to 
  refer to indefinite subjects, it is a common practice of natural 
  languages as well. As a matter of fact, the following sentence 
  illustrates a possible use of the word `substitute' in the Oxford 
  Advanced Learner's Dictionary of Current English ([Ho]): `Mr X  
  substituted for the teacher who was in hospital'. Of course, the 
  following example would have served us even better: `The teacher 
  substituted for Mr X who was in hospital'.}  Well, the substitution  
 is a formal, entirely language--independent operation, its result does 
not depend on the specific structure of the sentence which it is 
 applied to, the variable 
  guarantees 
that 
the subject of the new sentence we obtain after the substitution will 
indeed be the name of the given object.
Moreover, this transformation is also general in the sense that it, in 
principle, can be defined for formal and informal languages as well. 

\hspace{2mm}  We can now obtain a language--independent abstract form of the 
notion of self--reference  by changing the transformation in 
(3) to substitution in a suitable way (and using the fact that the 
quotation is only a particular way to give a name to a linguistic object): 

\vspace{1mm}

(6) \hfill   `applied to itself'  \ is equivalent to \ 
 `its name is substituted for the variable in it'.  \hfill \

\vspace{1mm}

\hspace{2mm} Accepting this interpretation, however, has 
its consequences.  Clearly, 
(6) implies directly that a linguistic entity, in order for it to be 
meaningfully applied to itself, has to be a sentence having a single 
variable in it. Therefore, to obtain the version of the Liar we seek,  
we should modify our original notion of heterologicality according to 
this new notion of self-reference. Well, we shall define {\it sentences 
with single variables} (not adjectives or adjectival phrases) to be 
heterological. It is fairly obvious that if we look for a sentence 
with a single variable that in a way expresses the same thing as an 
adjective $a$, then the most immediate candidate will be the sentence 
`$x$ is $a$'. Moreover,  (6) implies 
that we have to consider the 
sentence `$x$ is $a$' to be    
heterological just in the case that `\,``$x$ is $a$'' is $a$\,' 
is a false sentence. Since, 
with the notion of substitution at our disposal, we are already able 
to be more formal, let us give the exact definition:

\vspace{1mm}

\ \hfill 
\begin{minipage}[t]{13.2cm} A sentence with a single variable is 
called {\it heterological} if the new sentence \linebreak obtained 
by substituting the name of the sentence for the variable in 
it is false. 
\end{minipage} \hfill \

\vspace{1.5mm}

For example, since `$x$ is not a string of letters'  {\it is} 
a string of letters, 
the sentence `\,``$x$ is not a string of letters'' 
is not a string of letters\,' is false, consequently 
 `$x$ is not a string of letters' is heterological.    
 Now, the only thing 
that remains to do is to turn back to (2) and reformulate it using 
our modified notion. The sentence corresponding to (2) will be 

\vspace{1mm}

(7)  \hfill   `$x$ is heterological' is heterological.  \hfill \ 

\vspace{1mm}

\hspace{3mm}Using our definition of heterologicality above, 
`$x$ is heterological' takes the following form: 
 `the new sentence obtained by substituting 
     the name of the sentence 
     $x$ for the variable in it is false', 
which, in turn, by (7), provides the {\it Findlay paradox}\,: \footnote{J. Findlay used sentences of the same structure 
  to examine 
  \mbox{informally G\"{o}del's incompleteness theorem (cf. [F]).}}

\vspace{1mm}

\ \hfill 
\begin{minipage}[t]{11.3cm}  
 the new sentence obtained by substituting the name of the sentence
  `the new sentence obtained by 
substituting the name of the sentence
 $x$  for the variable in it is false' for the variable in it is false. 
\end{minipage} \hfill \

\vspace{2mm}

Note that just in the same way as we did in the case of the Quine 
paradox, it can directly be checked that this sentence indeed says of 
itself that it is false (and says nothing more) since it is built up in such a way that if we 
perform the substitution described in it, then we obtain the sentence 
itself, which is stated to be false. Moreover, it is a perfectly 
well formed sentence. Finally, its paradoxical character does not depend on any contingent linguistic fact, it is entirely 
language--independent so can immediately be translated to any language,
 be it natural or formal. Thus, at last,  having an impeccable 
 common language version of  
our original paradox, we can turn to embedding it to an abstract formal 
logico--linguistic system in order to examine the possible ways to 
resolve it.

\vspace{1mm}

\pagebreak

\begin{center}

 {\large\bf           2. The formal resolution of the paradoxes  }
\end{center} 

\vspace{0mm}

{\bf a. The abstract system of Smullyan}

\vspace{1mm}

\hspace{2mm}The Findlay paradox implicitly contains the minimal 
requirements an 
abstract system within which the paradox can be formulated should 
meet. In fact, it has to contain the formal versions of elements of 
natural languages used for setting down the paradox. These are the 
formal versions of any string of letters (they will be called 
{\it expressions}), those of the (meaningful) {\it sentences}, which 
can be true 
or false, and the sentences having single variables (in formal terms 
these two sets together are called {\it formulas}), a set of objects that 
will play the role of the {\it names} of the linguistic objects, and two 
mappings for representing the procedure of {\it naming} and that of the 
{\it substitution}. Finally, the notion of truth can be represented 
formally by an arbitrary subset of sentences. The elements of this set 
we shall 
consider to be true. Our definition below is a slightly modified version of the 
notion appearing in [S2] (p.5).\footnote{Smullyan, 
in excellent papers and books,  
examined 
 G\"{o}del's theorem in original ways both from the point of view of 
various abstract formal systems and within     
the framework of  
ingenious logical puzzles, which are very 
interesting and enlightening both to the expert and the interested 
layman (see e.g. [S1], [S2], [S3]).} 
First we set some notation. 

\vspace{1mm}

{\bf Notation }

\vspace{0.5mm}

For any sets  $X$, $Y$, and $Z \subseteq X$,  and for any function  $f : X  \longrightarrow  Y$\,,  

\vspace{0.5mm}

 \ \hfill $\displaystyle{  X \sim Y = \{x \in X : x \not \in Y\},  }$
 \ \  \ 
$ \displaystyle{   f^*Z = \{y \in Y\,:\, f(x) = y
  \mbox{\ for some\ }  x \in Z\}}$. \hfill \ 

If $f$ is one-one, the inverse of $f$ will be denoted by 
$\displaystyle{ f^{-1}}$.

\vspace{1mm}

{\bf Definition }

\vspace{0.5mm}

$ {\mathscr S} = \langle E, S, F, N, g , s \rangle $ is an 
{\it abstract formal system} if 

\vspace{0.4mm}

\ (i)\ \ \ \begin{minipage}[t]{15cm}
  $\emptyset \neq S \subseteq F \subseteq E$ and 
$N$ is an arbitrary set such that 
$F \cap N = \emptyset$\,. $E$ is the set of 
{\it expressions}, while 
  $S, F$, and $N$ are called the sets of 
  {\it formulas, (formal) sentences}, and {\it  names} respectively. The elements  of $F \sim S$ are the {\it proper formulas}. 
 We set\,   
$\sim A = F\sim A$\,\, for  any $ A \subseteq F$\, and\, 
$\sim X = N \sim X$\, for  any $ X \subseteq N$\,. 

\vspace{0.5mm}

\end{minipage}

\vspace{0mm}

\ (ii)
\ \ \begin{minipage}[t]{15cm}
    $g$  is a one-one mapping from  $F$  onto  $N$. $g$ is the 
         {\it naming function of}  ${\mathscr S}$. 
For any subset \linebreak of $F$, we shall denote the image 
of this subset under 
$g$ by the boldface version of the letter denoting the subset 
concerned, that is, e.g. if  $H \subseteq F$, then  $g^*H$  will be  
denoted by  ${\mathbf H}$. 

\end{minipage}

\vspace{0.5mm}

\ (iii)
\ \ \begin{minipage}[t]{15cm}
    $s$  is a mapping from  $F \times N$  into $S$  such that 
 \ $s(\sigma,n)= \sigma$\, for any $\sigma \in S,\, n \in N$. \
$s$ is the {\it substitution} in ${\mathscr S}$. We shall denote 
 $s(\varphi,n)$  by  $\varphi[n]$  for any  $\varphi 
\in F, \, n \in N$. 

\end{minipage}

\vspace{2mm}

 \hspace{2mm}Now, let us try to formulate the formal version of the Findlay 
paradox in ${\mathscr S}$. Let  $T \subseteq S$  be arbitrary. We shall  consider $T$ to be the 
set of true sentences of ${\mathscr S}$ and, by the same token, the sentences 
outside $T$ to be the false ones. Further, in the case of
 any common language sentence $s$ having a single variable 
$x$, we shall abbreviate the result 
of substituting a linguistic phrase  $q$  for the variable  
$x$  in $s$ by $s(q)$. 
Just as we did so far, we shall denote the 
name of any common language expression $e$ by `$e$'. 
 Moreover, let us denote the phrase 
\  `the new sentence obtained by substituting the 
       name of the sentence $x$ for the variable in it is false'
by  $p$. 
Using these notations, the Findlay paradox, 
abbreviated by $f$ from now on, \mbox{takes the following form:} 

\vspace{0.5mm} 

(8)  \hfill        $   f = p($`$p$'$)$.  \hfill \

\vspace{0.5mm} 

\hspace{2mm}
  First of all, note that only {\it names} of objects can be substituted for 
a variable in a linguistic phrase.  Thus $x$ stands for a 
 {\it name} in the phrase 

\vspace{1mm}

(9) \hfill 
         the new sentence obtained by substituting the name of the 
                sentence $x$ for the variable in it.
\hfill \ 

\vspace{1mm}

Therefore, in (9), the name of the sentence $x$ {\it is just} 
$x$ \footnote{Similarly as e.g. the name of the sentence 
`\,This is obvious.\,' 
is the following 
object: \ `\,This is obvious.\,'\,.},  
 so in its 
formal version again the name of the formula  $x$  is just  $x$, 
consequently the formula itself is $g^{-1}(x)$ and 

\vspace{0.5mm}

(10) \hfill    the formal version of the whole phrase (9) is 
               $g^{-1}(x)[x]$\,.\hfill \

\vspace{0.5mm}

Consequently, using the terminology provided by  ${\mathscr S}$,  
$p$  corresponds to the following sentence: 

\vspace{1mm}

(11) \hfill        $g^{-1}(x)[x] \not \in T$\,.  \hfill \

This sentence will be denoted by  $\overline{p}$.

\hspace{2mm}   It is important to note that the phrase (9) is undefined, 
it `depends 
on $x$', that is, it does not denote a unique sentence, the object it 
refers to may change with the name of the sentence that is substituted 
for the variable $x$. This property is, of course, inherited by its 
formal version $g^{-1}(x)[x]$ as well as by the phrases they are 
contained in, namely  $p$  and  $\overline{p}$, so  $\overline{p}$  
becomes a single sentence only as a result of substituting an  
$n \in N$ (a formal name of an expression) for the variable $x$ in it.  

\hspace{2mm}  Well, of course,  $\overline{p}$ is generally not an 
expression {\it belonging to} ${\mathscr S}$, 
but a statement {\it about} ${\mathscr S}$\,! It is still a natural language 
reformulation of  $p$  in terms of the abstract formal system  
${\mathscr S}$. Consequently, 
not belonging to the set of expressions of ${\mathscr S}$, it does not 
posses a formal name (it is not in the domain of the naming function $g$), so 
there is no formal version of \, `\,$\overline{p}$\,'. It follows,
 then, that we 
cannot continue the process of formalization, which, by (8), would be 
the substitution of this name for the variable  $x$  in  
$\overline{p}$   in order to obtain the formal version 
of  $f$. In other words, the paradox simply disappears since it cannot 
be formulated at all.\footnote{A. Tarski examined first systematically the possibility of 
eliminating the paradox along these lines (cf. [T]).} Consequently, we get 
to a crucial point in the reconstruction of the paradox within our 
formal framework.

\vspace{2mm}

{\bf b. The Formal Liar}

\vspace{1mm}

\hspace{2mm}   For obvious reasons, we do not want to breathe 
 life into 
our dead 
paradox (for that matter, even if we wanted to, we could not). What we do 
want is to reformulate it into a statement about the expressive power 
of our formal system, a statement to the effect that the absence of 
limits concerning this expressing power leads to a contradiction. 
Actually, the fact that prevents us from continuing our
 procedure can, 
from the opposite point of view, be taken as a condition under which 
this procedure can, in fact, be accomplished. Indeed, though $\overline{p}$ does 
not belong to ${\mathscr S}$, there might exist elements of ${\mathscr 
S}$ that can, in some way, 
represent it, that is, play the role of  $\overline{p}$   {\it within } 
${\mathscr S}$ ! Of course, 
these elements, in some sense, should `express the same state of 
affairs' as $\overline{p}$ . Now, recall that  
$\overline{p}$  is a sentence with an indefinite subject so it 
cannot be true or false since it `depends on  $x$'.  $\overline{p}$  
  becomes a 
single sentence only after substituting an  $n \in N$ for the variable  
$x$   
in it. Well, those formal expressions that have the corresponding 
property are just the proper formulas of ${\mathscr S}$. 
So we are seeking a 
formula, say  
$\pi$, `to the same effect' as  $\overline{p}$\,, which obviously means that  $\pi$ and  $\overline{p}$    
have to be true and false at the same time, in other words, they have 
to be true exactly for the same formulas, that is, by (11), 

\vspace{1mm}

(12) \hfill for any  $ n \in N$, \ \   $\pi[n] \in T$ \   iff  
\  $g^{-1}(n)[n] \not \in T$. \hfill \

\vspace{1mm}

Our arguments above made it clear that the existence of such a formula 
characterizes an abstract formal system in a fundamental way 
guaranteeing that a given object 
(the set of expressions whose  
 names
\footnote{Let us stress again: we use the   
 {\it names} of objects to talk about them.} 
 satisfy (11)\,) 
can be talked about within the system. Thus the expansion of this 
 property to any set of objects is related to the expressive 
 power (or, as it is often called, the strength) of a formal system.

\vspace{1mm}

{\bf Definition }

\vspace{0.5mm}

Let  ${\mathscr S}$  be an abstract formal system, and let 
$T \subseteq S$\,, $X 
\subseteq N$  be arbitrary. 
We say that  $X$  is   {\it $T$--representable} ({\it in} ${\mathscr S}$)  if  
there is a  $\varphi \in F$  such that for every  $n \in N$, 

\vspace{0.5mm}

\ \hfill    $\varphi[n] \in T$  \  iff  \   $n \in X$\,.  \hfill \ 

\vspace{0.7mm}

The formula  $\varphi$  is said to  {\it $T$--represent} $X$ 
({\it in} ${\mathscr S}$). 

\vspace{1mm}

Well, if   $\pi$  represents the set 
$\{n \in N : g^{-1}(n)[n] \not \in T \}$, that is, 
the condition (12) holds, then we can consider  $\pi$  as a representative 
of  $\overline{p}$  
 within  ${\mathscr S}$  and can continue the process of 
formalization using 
it as the formal version of $p$ (recall that $\overline{p}$ is the 
reformulation of $p$ 
in terms of  ${\mathscr S}$) to obtain the formal 
version of \  $f = p($`$p$'$)$. In fact, in this case, the formal object 
corresponding to  `$p$'  is  $g(\pi)$, so the one corresponding to  $f$ is  
$\pi[g(\pi)]$. This  is the sentence of  ${\mathscr S}$  corresponding to 
the Findlay 
paradox, which is, in turn, nothing else than a reformulation of 
the Liar. 
Consequently, we have obtained the 
formal version of the Liar. With the notations above, 

\vspace{0.5mm}

\ \hfill the  {\it Formal Liar} is the
 sentence $\lambda = \pi[g(\pi)]$.\hfill \

\vspace{0.5mm}

This is the formal 
sentence we sought. Like its informal counterpart, it turns 
 out to  
 be neither true nor false witnessing the fact 
that  
the formula $\pi$ it is built on cannot exist. 
The statement to this effect is 
the  resolution of our paradox in a formal setting.    
Note that the proof is, in fact, the formalization of the 
informal argument: the Formal Liar is true  iff   it is false. 

\vspace{1mm}

{\bf Proposition}  (The Liar Theorem)

\vspace{1mm}

{\it Let  ${\mathscr S}$  be an abstract formal system and $T \subseteq S$. 
The set \ $\{n \in N : g^{-1}(n)[n] \not \in T \}$ \ is not 
$T$--representable.}

\vspace{1mm} 

{\sc Proof.}

Let us suppose that, on the contrary, there is a $\pi \in F$ such that 
$\pi$  $T$--represents the set \linebreak $\{n \in N : g^{-1}(n)[n] \not \in T \}$, 
that is, 
$\pi[n] \in T$  \  iff \   $ g^{-1}(n)[n] \not \in T $.  Let  
$\lambda = \pi[g(\pi)]$  and  $n = g(\pi)$. Then  $\lambda = 
\pi[g(\pi)] \in T$  \  iff \  
$\pi[n] \in T$   \  iff  \   $g^{-1}(n)[n] \not \in T $  \  iff  \ 
$g^{-1}(g(\pi))[g(\pi)] \not \in T$ \  iff \  $\lambda = \pi[g(\pi)] 
\not \in T$, which is a contradiction. 

\vspace{1mm}

\hspace{2mm}   So we have the formal version of the Liar at last. To be 
entirely 
sincere, the paradox in this form is not too exciting since it is not 
 easy to interpret in informal terms, or, simply, to understand 
what it is all about. In order to grasp its essence, we should make 
it more transparent by giving a name 
to the fundamental 
notion lying at the heart of the whole matter we examine. 

\vspace{3mm}

{\bf c. The Generalized Liar Theorem}

\vspace{0.5mm}

\hspace{2mm}   Before analyzing the Liar theorem in order to understand 
it better, 
however, we should make a remark. Actually, there is an important 
point to stress here. When we formulated the definition of 
representability, we have taken for granted that the answer to the 
question whether a given statement is satisfactory, agreeable, 
acceptable etc. or not depends on its  {\it truth}. The truth of a 
statement, on the other hand, means some kind of proper correspondence 
or adequacy between the state of affairs and the statement concerned.  So 
far so good. But let us stop here for a moment. Truth is, by no means, 
the only concept connecting facts to statements in some way or other, 
thus it is not the only possibility to chose from for formulating the 
definition of representability. Actually, there could exist and do 
exist other `measures of adequacy' than truth that lend themselves to 
formalization more or less readily. It is enough to think of the 
notions of probability or confirmability (to different extents). These 
are, clearly, `weaker' notions than truth. In the other direction, we 
can find the `stronger' ones the most obvious choice of which is, of 
course, that of provability. To be more precise, it would be better 
to speak 
about  provability in plural since one can introduce, in a very 
natural way, a couple of different notions of provability widening or 
narrowing the circle of stipulations as to which types of logical 
derivations would be accepted as a proof. Actually, we shall briefly 
touch upon one of these possible modifications of the generally 
accepted notion of provability below.  

\hspace{2mm}   Anyway, it is natural to generalize our treatment to 
include any 
one of these possible `measures of adequacy' and not to consider our 
definition above as the only possible definition of representability,  
but, rather, as one of them. Accordingly, in the same way as we did in 
the case of true sentences, we shall represent formally the sets of 
sentences satisfying different possible adequacy conditions simply by 
arbitrary subsets of the set of sentences 
 and, 
until further 
notice, shall  not associate any special meaning with them. 
As a reminder of this fact, we repeat our representability 
definition using a letter in it without any connotation\,: 

\vspace{1mm}

{\bf Definition}

\vspace{0.5mm}

Let  ${\mathscr S}$  be an abstract formal system, and let  
$A \subseteq S$,  $X \subseteq N$  be arbitrary. We say that  
$X$  is  $A$--{\it representable} ({\it in} ${\mathscr S}$)  
if  there is 
a  $\varphi \in F$  such that  for every  $n 
\in N$, 

\vspace{0.5mm}

\ \hfill $n \in X$  \  iff \  $\varphi[n] \in A$\,.  \hfill \ 

\vspace{0.5mm}

The formula $ \varphi$  is said to $A$--{\it represent} $X$
 ({\it in} ${\mathscr S}$). 

\vspace{1mm}  

\hspace{2mm}   Now, we can turn back to the Liar theorem. 
It is, no doubt, the 
expression $g^{-1}(n)[n]$  which prevents us from interpreting it 
informally in an easily comprehensible way. Everything would seem 
much simpler if we were able to get rid of it. Certainly, 
we cannot avoid handling it in some way or other. 
 For recall  
that, by (6), (9), and (10), $g^{-1}(n)[n]$  is just the formal 
version of the phrase 
     `the new sentence obtained by applying 
the sentence $n$ to itself', 
 that is,  $ g^{-1}(n)[n]$   is nothing 
else than the 
formalization of  {\it self--reference} and, clearly, due 
to its central role in the paradox, it will appear 
in any formulation of the Liar
 blurring the resulting picture unless we 
restrict ourselves to formal systems that can circumvent it.    
 Well, in those   
formal systems which can do without self--reference, it 
is possible to talk about everything without explicit 
self--reference that can be talked 
 about at all, 
that is,  
in our 
formal terms, any set of names 
 defined by a 
representable set through self-reference is representable 
itself.\footnote{ This 
property corresponds to the condition  
formulated in the {\it fixed--point} (or {\it diagonal}\,)\, lemma 
of mathematical logic. The most obvious examples of 
formal systems 
with this property are those induced by natural languages, 
which are  
rich enough to contain the expression corresponding to 
the general self--reference 
itself (cf. (9)). 
What is more important is that the same is true for 
the usual formal systems of 
arithmetic.}

\vspace{0.5mm}

{\bf Definition }

\vspace{0mm}

Let  ${\mathscr S}$  be an abstract formal system and  
$A \subseteq S$. 
${\mathscr S}$  is 
     {\it self-referential with respect to} $A$  if  for any 
   $X \subseteq N$,  

\vspace{-0.7mm}

   \ \ \hfill $\{n \in N : g(g^{-1}(n)[n]) \in X \}$  is $A$--representable whenever 
   $X$ is $A$--representable. \hfill \

\vspace{1mm}

\hspace{3mm}With the notion of self--referentiality at our disposal, 
we get to the last step of our  
formalization process. Recognizing that the 
set $T$ plays a double role in the Liar theorem, it is very natural to 
modify the wording and proof of this statement   
 in the only obvious way to obtain a theorem  
 and its proof about the relation between {\it two} sets of 
sentences in place 
of a less general proposition being about only a single one.  
This modification of the Liar theorem,  
together with the 
application of self--referentiality, yields the 
generalized formal version of the original Liar paradox.  

\vspace{0.5mm}

{\bf Theorem}  (The Generalized  Liar Theorem)

\vspace{0.5mm}

{\it Let  ${\mathscr S}$  be an abstract formal system, 
 $A \subseteq F\,,\, B \subseteq S$. 
 Let  ${\mathscr S}$  be self--referential with respect to  
$B$  and suppose that  $\sim {\mathbf A}$  is  $B$--representable. 
Then \  $ S \cap A \neq B $\,. } 

\vspace{0.5mm}

{\sc Proof.}

Since  $\sim \!{\mathbf A}$  is  $B$--representable 
and  ${\mathscr S}$  is self--referential 
with respect to $B$, \mbox{$\{n  \in \! N \!: g(g^{-1}(n)[n]) \! 
\not \in \!
{\mathbf A} \}$}
is also $B$--representable, that is, there is a  $\pi \in F$  such that 
for any  $n \in N$,  $\pi[n] \in B$   iff   
$g(g^{-1}(n)[n]) \not \in {\mathbf A}$. Let  $n = g(\pi)$  and 
 $\lambda = 
\pi[g(\pi)]$. Then, 
on the one hand, $\lambda \in S$. On the other hand,  $\lambda = 
\pi[g(\pi)] \in 
B$   iff   $\pi[n] \in B$   iff   $g(g^{-1}(n)[n]) \not \in {\mathbf A}$
  iff   
$g^{-1}(n)[n] \not \in A$   iff   $g^{-1}(g(\pi))[g(\pi)] \not \in A$ 
  iff  
$\lambda = \pi[g(\pi)] \not \in A$.  That is, we have a 
$\lambda \in S$, such that 
    $\lambda \in B$ \   iff  \  $\lambda \not \in A$\,,  
which was to be proved.

\vspace{0mm}

\hspace{2mm}   Of course, in such an abstract wording, this theorem does 
not seem 
to tell much about the relevant formal systems of mathematics. Before 
applying it to the logical notions that are the usual objects of 
metamathematical investigations, however, we would like to illustrate 
the way it can  be used by applying it to 
the notion of provability 
{\it in less than a given 
number of steps} not usually examined despite its very interesting 
intuitive meaning. 

\hspace{3mm}Mimicking the informal systems of 
mathematics, we can define the set of provable sentences in an 
abstract formal system ${\mathscr P}$
 along the following lines. After selecting a 
set of sentences that are called axioms and a set of rules for 
inferring an expression from a finite set of other ones, a  {\it 
proof} is 
defined as a finite sequence of expressions every element of which is an 
axiom or can be inferred from the set of elements preceding the given 
one. A sentence is defined to be  {\it provable}, in other words to be a 
 {\it theorem}, if it is the last element of some proof. Now, 
for any positive integer  $n$, let us define 
the set  $P_n$  as that of the sentences provable in less than  $n$  steps, 
that is, appearing as the last element of a proof consisting of less 
than  $n$  expressions. Further, let  $P$ be the set of all provable 
sentences of  ${\mathscr P}$. Clearly, $P_n \subseteq P$  for any  $n$. 
Thus if we interpret  $A$  
and  $B$  in the Generalized Liar theorem   as  
$P_n$  and  $P$  respectively, 
then we have 

\vspace{1mm}

{\bf Proposition }

\vspace{1mm}

Let ${\mathscr P}$ be  self--referential with respect to  $P$. 
For any positive integer $n$, 
 if  $\sim {\mathbf P}_n$  is  
$P$--representable, 
then  

\vspace{-3mm}

\ \hfill 
    $P \sim P_n \neq \emptyset$. 
\hfill \ 

\vspace{2mm}

\hspace{2mm}The proof of the fact that in the most important cases of the systems 
of formalized arithmetic the conditions of the proposition actually hold 
for any positive integer is essentially simple. Yet, being a bit too 
technical, it does not fit into our non--technical exposition. 
Nevertheless, the conclusion of the proposition is worth spelling out in a 
little more detail since its informal interpretation reveals an 
interesting feature of the system considered. Actually, it means that 
for any positive integer  $n$, however large it is, there are theorems of the 
system that can only be proved in more than  $n$  steps, that is, there 
are arbitrary long proofs establishing new theorems. This, on the one 
hand, means some kind of practical incompleteness since there are 
theorems of the system whose demonstration needs an astronomical 
period of time, consequently, although they are theoretically 
provable, will never be proved for banal practical reasons: 
mankind simply has not and will never have enough time to prove 
them. From another point of view, this very property corresponds 
formally to the `unboundedness' of the formal system considered, 
that is, to the fact that there are new theorems being arbitrarily 
`far'\,\footnote{As a matter of fact, it is a little more than a 
suggestive metaphor to say that two sentences are  {\it far} from 
each other 
since the function assigning to any two sentences the length of the 
proof deriving the second sentence from the first one is a possible 
generalization of the notion of distance. In fact, it essentially 
coincides with a notion of functional analysis, namely that of the 
 {\it quasi--pseudometric}, which is a metric 
without the requirement of  
symmetry satisfying, among others, the triangle inequality.} from the 
axioms thus asserts formally its inexhaustibility, complexity, or 
richness.

\hspace{2mm}   Now, let us turn to those special abstract formal 
systems, that constitute the conceptual framework for investigating 
the theoretical aspects of mathematical activity.

\vspace{-1mm}

\begin{center} 
 {\large\bf   3.  Limitations of logical systems}
\end{center}

\vspace{-1mm}

\hspace{2mm}   In order to apply the  Generalized Liar theorem to the real 
systems of mathematics, we should formulate the minimum requirements 
any ordinary  
 logic 
has to comply with. Obviously, such a system ought to have two 
distinguished subsets of sentences that can be interpreted as those of 
 {\it provable} and  {\it true} sentences respectively and the system has 
to be able 
to express the basic logical operation of  {\it negation}. We shall explain informally the basic properties characterizing the different logical systems after defining them 
formally.  

\vspace{1mm}

{\bf Definition }

\vspace{1mm}

(i)  \begin{minipage}[t]{15cm} We say that an abstract formal system 
  ${\mathscr S}$ is a {\it logical system} ({\it with respect to}  $P$
   {\it and}  $T$) if  $P, T \subseteq S$
and, for any  
  $\varphi \in F$, there is a $\varphi' \in F$, called 
the  {\it negation of} $\varphi$  (or  {\it  not}  $\varphi$) such that, 
 for any $n \in N$, 

\vspace{0.8mm}

    (a) \  $\varphi[n] \in T$ \  iff  \
 $\varphi'[n] \not \in T$\, \ \ \  \ 
  (b)  \  $\varphi[n]  \in P$ 
\ iff \  $\varphi''[n] 
 \in P$ \hfill \

\vspace{1.2mm}

 (In particular, for any $\sigma \in S$, \ 
 $\sigma \in T$  iff  
 $\sigma' \not \in T$\,, \, $\sigma \in P$  iff  
 $\sigma'' \in P $.) 
 $P$   and  $T$ are called  the sets of  {\it provable} and  
{\it true} sentences of ${\mathscr S}$ respectively. 
We shall use the following notation for any set  $H \subseteq F$: \
$    H' = \{\varphi  : \varphi' \in H \}\,.$ 

\end{minipage} 

\vspace{1.5mm}

(ii) \hspace{0.2mm}A logical system   is   {\it consistent} if  
$ P \cap P' = \emptyset\,$, 
otherwise it is  {\it inconsistent}.

\vspace{0.5mm}

(iii) A logical system  is  {\it complete} if 
\ $  P \cup P' = S$, 
otherwise it is  {\it incomplete}. 

\vspace{0.5mm}

(iv) A logical system  is 
  {\it sound} if \   $P \subseteq T$. 

\vspace{0.5mm}

\hspace{2mm}Obviously, consistency is one of those conditions that  any meaningful 
 logical system should satisfy since this means that it is 
 free of contradiction: for  
 any sentence and its negation, {\it at most} one of them 
 belongs to the set of 
 provable sentences. Moreover, clearly, the set of provable sentences has to be chosen in such a way that the system remain sound; only the 
true sentences are permitted to be provable.  On the other hand, our main concern is  to 
examine the  completeness  of formal systems, the property 
which is, in some sense, a maximal requirement,  a dual of the consistency. In complete  
systems, 
for an arbitrary sentence and its negation,  {\it at least} one of them does belong to 
the set 
of provable sentences. The  importance of this notion stems from the fact that, since all the interesting formal systems  are sound, any one of them 
is complete iff $ P = T$\,,\,\footnote{Indeed, supposing 
 the soundness,   $ P \neq T $ implies that      
there is a $\sigma \not \in P$ such that $\sigma \in T$  \, (iff \ 
 $\sigma' \not \in T$\,), that is,  
$\sigma' \not \in P$. On the other hand, by $T \cup T' = S$\,,
 $\sigma \not \in P$ and    
$\sigma' \not \in P$ implies that 
 $\sigma \in T \sim P$ or 
$\sigma' \in T \sim P$ i.e. $P \neq T$\,.} that is,   
{\it all} true sentences are, in fact, provable.\footnote{
 Moreover the notion of completeness is  important from  a `practical' 
point of view  
as well. Indeed, in the case of the formal systems
 describing 
mathematical structures, the complete systems 
are decidable in the sense that there is an entirely mechanical rule to answer the question whether a sentence is a 
theorem or not. 
Actually, in  
these systems, the set of the finite sequences of expressions
 (which 
themselves are finite sequences of symbols chosen from a finite 
vocabulary) can be enumerated, and, as they mirror the world of informal 
mathematics, the notion of the proof in 
these
systems is defined in such a way that the fact whether a given 
sequence of expressions 
is a proof or not can be decided in a finite 
number of steps. So we simply enumerate all the proofs until 
our sentence or its negation appears as a last element of some 
proof since we know that one of them, being provable, must emerge 
sooner or later.} 
  
\hspace{2mm}Well, since $P \,, P' \subseteq S$\,,
$P \cup P' = S$  \ iff \  
$ (S \sim P) \cap (S \sim P') = 
\emptyset\,$. Thus   
a logical system  ${\mathscr S}$  is  complete  {\it and} consistent  
iff $(P \cap P') 
\cup \left(\rule[0.0mm]{0mm}{3.5mm}(S \sim P) \cap (S \sim P')\right) 
= \emptyset$   
iff   $S \sim P = P'$. 
Clearly, the relation between the last equation and the 
conclusion of the Generalized Liar theorem needs no comment, 
they together yield an abstract version 
of G\"{o}del's incompleteness theorem to the 
effect that, under suitable conditions concerning self-referentiality 
and representability, {\it consistency and completeness exclude each other.}
  Actually, restricting ourselves 
to logical systems,  
the Generalized Liar theorem 
takes a form that yields the abstract versions of three 
basic limitation theorems  
 of mathematical logic 
 demonstrating that they are {\it all}  
different manifestations of the {\it same} logical 
principle.\footnote{Smullyan investigated a wide variety of 
abstract limitation theorems in [S2]. It is interesting, however, 
that he has apparently not recognized 
their common logical structure or, at least, has not 
found it to be worth examining.}

\vspace{1mm}

{\bf Theorem}  (Abstract versions of theorems of G\"{o}del, Tarski, 
and Church)

\vspace{1mm}

{\it Let  ${\mathscr S}$  be an arbitrary logical system. 

\vspace{1mm}

{\rm (i)} {\rm (a)} Let us suppose that  ${\mathscr S}$  
is self-referential with 
respect to  $P$ and ${\mathbf P}$  is  
$P$--representable. 

\ \hfill 
If  ${\mathscr S}$  is consistent, then  ${\mathscr S}$  
is incomplete.  \hfill \

\vspace{0.5mm}

\hspace{5mm}{\rm (b)} Let us suppose that  ${\mathscr S}$  
is self-referential with 
respect to  T and  ${\mathbf P}$  is  
$T$--representable. 

\vspace{0.5mm}

\ \hfill If  ${\mathscr S}$  is sound, then  ${\mathscr S}$  is incomplete.  \hfill \

\vspace{1mm}

{\rm (ii)}\, Let us suppose that  ${\mathscr S}$  is self-referential with 
respect to  T. Then 

\vspace{0.5mm}

  \ \hfill  ${\mathbf T}$   is not $T$--representable. \hfill \

{\rm (iii)} Let us suppose that ${\mathscr S}$  is  self-referential with 
respect to  $P$. Then  

\vspace{0.5mm}

      \  \hfill $\sim {\mathbf P}$   is not
 $P$--representable.}  \hfill \ 

\vspace{1mm}

{\sc Proof.}

We list the substitutions needed to obtain the different items from 
the Generalized Liar theorem:

\vspace{1mm}

(i)(a):    $A =  \sim P$\,,  $B = P'$ \
(i)(b):    $A =  P$\,,  $B = T$   \ 
(ii):   $A =  T$,  $B = T$ \ \ 
(iii):   $A = P$\,,  $B = P$\,. 

\vspace{1mm}

Going into a little more detail, it directly follows from the definitions that, generally,
 for any set  $C \subseteq S$, if  ${\mathbf C}$ 
 is $T$--represented by  
$\varphi \in F$, then  $\sim{\mathbf C}$   is $T$--represented by  
$\varphi' \in F$ and  if  ${\mathbf C}$  is $P$--represented 
by  $\varphi \in F$, then  ${\mathbf C}$  is $P'$--represented by  
$\varphi' \in F$\,. Further, self--referentiality of 
${\mathscr S}$ with respect to $P$ implies its self--referentiality 
 with respect to $P'$. Using these facts where needed, we have:

\vspace{0mm}

(i)(a) $\sim \sim {\mathbf P} =  {\mathbf P}$\,, so, by the 
Generalized Liar, $S \sim P \neq  P'\,$, which in turn, by our remarks preceding the theorem, means that    
 consistency implies incompleteness. 

(i)(b) By the Generalized Liar,  
$P \neq T$\,. 
${\mathscr S}$  is sound, thus ${\mathscr S}$ is incomplete (cf. 
footnote 22). 

(ii)  By the Generalized Liar, $T$--representability 
of $\sim \!{\mathbf T}$ implies  \mbox{$T \neq  T$.} 

\vspace{0mm}

(iii) By the Generalized Liar,  $P$--representability 
of  $\!$ $\!$ $\sim \!{\mathbf P}$ implies  \mbox{$P \neq P$.} 

\vspace{0.5mm}

\hspace{2mm}   Let us underline that, 
 the theorem above is, in a sense, exhaustive, it  
contains the generalizations of all the known 
main limitation theorems on truth and provability 
that can be formulated 
on the level of abstraction implicitly 
defined by the formalization of the ordinary--language Liar.\footnote{Indeed, the only missing important limitation 
theorem, perhaps the most important one, 
G\" odel's theorem on the unprovability of consistency 
(to the effect that `arithmetic cannot prove 
its own consistency') cannot   
be taken as a different limitation theorem on the level 
of our formalization since it is, essentially, 
the  formal version  of 
 the  incompleteness theorem itself given on a still deeper level, 
 {\it within} formal arithmetic.} Actually, 
in addition to its first main item consisting of 
two abstract versions of G\"odel's incompleteness  
theorem describing the relation between 
{\it  provability and refutability}
\footnote{A sentence is refutable if its negation is provable. 
Note that this statement 
 is a purely syntactical 
one having essentially nothing to do with the notion of truth.}  
and that  
 between {\it  provability and truth}\,, respectively, it contains, 
 as its second and third main 
items, a generalization of the Tarski theorem on the {\it undefinability of 
 truth}\,,\footnote{Informally this result can be interpreted as stating that \mbox{the notion of being true in  
${\mathscr S}$
 cannot$\!$ be$\!$ defined within $\!{\mathscr S}\!$.}} and 
 an abstract 
variant of the theorem of Church on the 
{\it undecidability of  provability}
$\!\!$
\footnote{In the 
case of the 
usual formal systems describing the arithmetic of natural numbers, the informal content of this theorem is that there is 
no mechanical way to decide whether a given sentence is a 
theorem or not. In fact, in the case of these systems (in which $N$ is 
a subset of natural numbers), there are 
several equivalent formal notions devised to formalize 
the informal notion of decidability of 
the question 
whether a given natural number  belongs to a given set of natural 
numbers. One of them is the property 
 that {\it both} the set under consideration and its 
complement are $P$--representable. An argument completely  
analogous to 
that given in footnote 23 may shed some light on this 
  connection between decidability and provability,  
at least in one direction.}.\footnote{The original theorems can be found 
in every textbook on mathematical logic, see e.g. [E] or [M].} 
Moreover, 
covering every sensible cast of roles for the Generalized 
Liar theorem,  
its  
proof indicates that the theorems of G\" odel, Tarski, and 
Church are just the only possible relevant 
limitation results 
given in terms of truth and provability alone 
 that can be considered as some direct reformulations 
of the Liar paradox within the conceptual framework of 
modern mathematical logic. 

\hspace{2mm}   As a matter of fact, on the other hand, this theorem constitutes the 
final point 
where we can get to solely on the basis of the Liar paradox without 
entering into the detailed analysis of individual formal systems. The 
theorem above describes the general  {\it conceptual structure} of the 
limitation theorems concerning the formal systems of mathematics,  
making explicit the different roles played, so to speak, by some 
 {\it general} laws of logic (represented by the Liar) on the one hand, and, 
on the other, by the special features of the particular
 system concerned. 
Actually, the theorem shows that, in the case of any given system, its 
limitations depend on its expressive power regarding 
self-referentiality and representability. In other words, examining 
 completeness, decidability, or the capability to express the notion 
of truth   
in the case 
of any segment of  real 
mathematical activity (as e.g. the arithmetic of natural numbers), 
where the expressions defined by the structure under investigation are 
given, the main question reduces to the more or less technical one as 
to whether we can find a naming function such that the resulting logical 
system satisfies the conditions of the theorem concerning 
self-referentiality and representability. From this general point of 
view, the essence of G\"{o}del's original result on the incompleteness of 
formal arithmetic consists in demonstrating that such a naming 
function, in fact, exists in this special case,   
so that the theorem can indeed be applied to the logical systems 
of 
arithmetic.\footnote{G\"{o}del discovered the ingenious method of associating 
   {\it numbers} with expressions in such a way that the resulting naming 
  function has the most important property of the quotation as the 
  naming of ordinary languages, viz. that the name of an expression 
  contains every important information of the named object (clearly, 
  the quoted version of a common language expression is ideal in this respect) thus the logical system endowed with such a naming function inherits 
the expressive power of natural languages, among others, their 
self--referentiality and the fact that $P$ is representable.  
Loosely speaking, 
the expressive power 
  of arithmetic can be compared to that of the natural languages 
  since it can also be forced to talk about itself.  
   The discovery of this revolutionary 
  technique is often  likened to the invention of Cartesian coordinate 
  geometry, that is, `G\"{o}del invented what might be called  
  {\it co-ordinate metamathematics}' (R.B. Braithwaite). 
This method of 
so--called G\"{o}del  numbering can easily be 
modified to provide similar results in the case of 
important mathematical theories other than arithmetic. Detailed 
description of 
G\" odel numbering and the proof of the representability 
of $P$ in the usual systems of 
arithmetic 
can be found 
 in every treatise 
on mathematic  
  logic, among others e.g. in [E] and [M].  
On the other hand, although the notion of self--referentiality  
 is not 
usually treated in the textbooks explicitly, 
to check that it holds for the systems considered is only a 
matter of a simple and straightforward calculation.}

\hspace{2mm} Finally, let us apply our theorem 
  in a nonstandard way in the 
opposite 
direction. Presburger has shown that the additive number 
theory  
(the weakened version of the standard Peano arithmetic obtained 
by the omission of  multiplication)   
is complete \mbox{(cf. [C] p. 43.)}  In 
this system the set of all expressions can be enumerated, and
 the question 
 whether a given element of the resulting sequence is 
 a provable sentence or not 
 can be decided (cf. footnote 23). 
So the naming function  $g$  can be defined by recursion 
in the following way. For any expression $e$, we  
set $g(e)$, according to whether it is 
a provable sentence or not,  
to be the smallest even or odd natural  
number, respectively, that has not already been chosen.  
Now, $P = T$\,, since the system considered is  
complete. On the other hand, both the 
set of all even natural numbers and that of all 
 odd ones 
{\it are} representable\,\nolinebreak
\footnote{In fact, the system contains  
the formulas $(\exists y)(y + y = x)$\, and 
$(\forall y)(y + y \neq x)$.}.
Consequently, 
in this system  ${\mathbf T}$,  ${\mathbf P}$,  and  $\sim {\mathbf P}$
 are both  $P$--  and $T$--representable 
(so, among others, the truth {\it is} definable\,\footnote 
{That is,  we have a very natural 
example of a system which `can define 
its own truth'; actually, this is a more transparent and 
self-explanatory example than  Myhill's one, cf.[My].}) and, 
therefore, by the 
theorem, this system is self--referential with respect to neither 
$T$ nor $P$\,. 

\hspace{2mm}   In the light of this example and the fact that 
G\"{o}del's 
incompleteness result provides examples of self-referential 
logical systems in which $P$ is  
representable (thus the system is incomplete), it would be 
interesting to find self-referential complete logical systems  
(if they exist at all), which will yield a complete description of 
the mutual correspondence between the general properties of 
self-referentiality, representability and completeness.  

\vspace{0.5mm}

{\bf Acknowledgement} I am grateful for  support provided 
\mbox{by Hungarian NSF grant No. T16448} and T023234.

\vspace{0.5mm}

{\bf References}

[A] J.L. Austin, {\it How to Do Thing with Words}, 
Oxford University Press, 
Oxford, 1980. 

[C] C.C. Chang, H.J. Keisler, {\it Model Theory}, North Holland 
Publ. Comp., Amsterdam, 1973.

[E] H.B. Enderton, {\it A Mathematical Introduction to Logic}, Academic 
Press, New York, 1972. 

[F] J.N. Findlay, {\it Goedelian Sentences: A Non-numerical Approach}, Mind, 
Vol.51, pp. 259-65.

[G] K. G\"{o}del, {\it On Formally Undecidable Propositions of Principia 
Mathematica and Related systems}, trans. by B. Meltzer, with an 
intr. \mbox{by R.B. Braithwaite, Oliver and Boyd, Edinburgh, 1962.} 

[H] D.F. Hofstadter, {\it Metamagical Themas}, Basic Books, Inc., 
Publishers, New York, 1985.

[Ho] A.S. Hornby, {\it Oxford Advanced Learners's Dictionary of Current 
English}, Oxford Univ. Press, 1977.

[K] A. Koestler, {\it The Act of Creation}, 
Dell Publishing Co. Inc., New
York, 1967.

[M] E. Mendelson, {\it Introduction to Mathematical Logic}, 
D. Van Nostrand Company, Inc., Princeton, 1965.  

[My] J. Myhill, {\it A System Which Can Define its Own Truth}, 
Fundamenta Mathematicae, \linebreak Vol. 37, 1950, pp. 190--192.

\mbox{[Q] W.V. Quine, {\it The Ways of Paradox and Other Essays}, Harvard 
Univ. Press, Cambridge, 1979.} 

[Q1] W.V. Quine, {\it Concatenation as a Basis for Arithmetic}, 
The Journal of Symbolic Logic, Volume 11, 1946. 
pp. 55-67. 

[R] R. Rucker, {\it Infinity and the Mind}, Bantam Books, Toronto, 1983.

[R] G. Ryle, {\it Heterologicality, Analysis}, Vol.11, No.3, 1951. 

[S1] R.M. Smullyan, {\it Languages in which Self Reference is Possible}, 
The Journal of Symbolic Logic, Volume 22, Number 1, March 1957, 
pp. 55-67. 

\mbox{[S2] R.M. Smullyan, {\it G\"{o}del's Incompleteness Theorems}, Oxford Univ. 
Press, New York, 1992.}

\mbox{[S3] R.M. Smullyan, {\it What is the Name of this Book?}, Prentice 
Hall, Englewood Cliffs, N.J., 1978.} 

[T] A. Tarski, {\it The Semantic Conception of Truth},
\mbox{ Philosoph. and Phenomenolog. Research
4 (1944).} 

\vspace{2mm}

{\footnotesize Technical University of Budapest

Department of Algebra

1111 Stoczek u. 2. H \'{e}p. 5.em.

Budapest, Hungary  

\vspace{-0.6mm}

e-mail: sereny@math.bme.hu

\end{document}